# DESCRIPTION OF THE SET OF ADMISSIBLE PIECEWISE-LINEAR ROUTES WITH $n$ TURNS IN THE THREE-DIMENSIONAL CASE

**V. N. Nefedov**

**Absrtract**
We consider piecewise-linear polygonal chains connecting two given points $A$, $B \in \mathbb{R}^3$ and consisting of exactly $n+1$ segments (i.e., having $n$ turning points). The absolute value of the turning angle at each interior point is bounded by a given number $\varphi \in (0, \pi)$. Under the condition $n\varphi \leq \pi$, we describe the set to which all interior vertices of such a polygonal chain belong (Theorem 1). It is proved that for any point $B^{(1)}$ from this set, there exists a polygonal chain with the specified parameters (Lemma 1). Based on these results, we obtain an explicit formula describing the set of all admissible sequences $\left(B^{(1)}, \dots, B^{(n)}\right)$ of the angular points of the polygonal chain. The obtained description can serve as a basis for constructing algorithms to enumerate admissible polygonal chains and to solve optimization problems for an objective function that accounts for the cost of traversing the segments and the cost of turns.

## Introduction

There are practically significant problems in which it is required to connect two given points in three-dimensional space by a piecewise-linear polygonal chain. In this case, there are constraints on the number of segments (exactly $n$) and on the absolute value of the turning angles at the breakpoints. The goal is to describe the set of admissible polygonal chains for their subsequent enumeration (in the discrete case) or approximation (in the continuous case). Based on such a description, one can solve optimization problems for an objective function that accounts for the cost of traversing the segments and the cost of turns.
In this work, we prove two main statements – Theorem 1 and Lemma 1, which describe the set of admissible polygonal chains. On the basis of these results, we derive a formula for the set of all possible sequences of angular points (Theorem 2). Analogous results for the two-dimensional case were previously presented in [1].

## 1. Basic definitions

### 1.1. The angle between vectors in three-dimensional space

In this work, we will not be concerned with the direction of the angle, as was the case in the planar two-dimensional scenario (see [1, 2]). That is, we assume that

$$\forall A, B \in \mathbb{R}^3 \quad \angle(A, B) \geq 0.$$

The absolute value of the angle $\angle(A, B)$ is uniquely determined by the value of the cosine of this angle (due to the strict monotonicity of the function $y = \cos x$ on $[0, \pi]$). As a result, in the case $A \neq \mathbf{0}, B \neq \mathbf{0}$, where $\mathbf{0} = (0,0,0) \in \mathbb{R}^3$, we obtain the formula

$$\angle(A, B) = \arccos \frac{(A,B)}{|A||B|} \in [0, \pi],$$

from which, for example, it follows that $\forall t_1, t_2 > 0$:

$$\angle(A, B) = \angle(t_1 B, t_2 A), \angle(A, B) = \angle(t_1 A, t_2 B) = \angle(-t_1 A, -t_2 B). \quad (1)$$

For zero vectors, to ensure greater generality in reasoning, we assume that

$$\forall A, B \in \mathbb{R}^3 \quad \angle(A, \mathbf{0}) = \angle(\mathbf{0}, A) = 0.$$

In this case, equalities (1) also hold when $A = \mathbf{0}$ or $B = \mathbf{0}$.

Using the properties of angles in Euclidean space and the triangle inequality for vector lengths, it is not difficult to show that the following triangle inequality for angles holds [3].

**Statement 1.** *For any vectors $A, B, C \in \mathbb{R}^3$, in the case $B \neq \mathbf{0}$, the inequality holds:*

$$\angle(A, C) \leq \angle(A, B) + \angle(B, C).$$

### 1.2. Angular sector and spherical torus in three-dimensional space

Let $A, B \in \mathbb{R}^3$, $A \neq B$, $\varphi \in (0, \ \pi]$. Let us introduce the set (spherical torus)

$$\mathbf{S}(A, B, \varphi) = \{C \in \mathbb{R}^3 \backslash \{A, B\} | \angle(B - C, C - A) \in [0, \ \varphi)\}. \quad (2)$$

To understand the structure of this set, consider, for an arbitrary point $D \in \mathbb{R}^3 \backslash AB$, where $AB = \{A + t(B - A) | t \in \mathbb{R}\}$ (with $AB$ being the line passing through points $A$ and $B$), the set $\mathbf{\Pi}(A, B, D)$ – the plane passing through points $A, B, D$, i.e. [4]:

$$\mathbf{\Pi}(A, B, D) = \text{aff}\,\{A, B, D\} = \{t_1 A + t_2 B + t_3 D | t_1, t_2, t_3 \in \mathbb{R}, t_1 + t_2 + t_3 = 1\}.$$

The line $AB$ divides the plane $\mathbf{\Pi}(A, B, D)$ into two half-planes: $\breve{\mathbf{\Pi}}(A, B, D)$ – conventionally the right one, and $\widehat{\mathbf{\Pi}}(A, B, D)$ – conventionally the left one, which share a common boundary – the line $AB$. For example, if $E \in \mathbb{R}^3$ is a non-zero vector perpendicular to the plane $\mathbf{\Pi}(A, B, D)$, then we can take $\breve{\mathbf{\Pi}}(A, B, D)$ as the set of all points $C \in \mathbf{\Pi}(A, B, D)$ such that either $C \in AB$, or the vectors $B - C, C - A$ and $E$ form a right-handed triple of vectors, i.e.:

$$\begin{vmatrix} B - C \\ C - A \\ E \end{vmatrix} > 0.$$

This leads to the following formulas:

$$\breve{\mathbf{\Pi}}(A, B, D) = \left\{ C \in \mathbf{\Pi}(A, B, D) \middle| \begin{vmatrix} B - C \\ C - A \\ E \end{vmatrix} \geq 0 \right\},$$

$$\hat{\mathbf{\Pi}}(A,B,D)=\left\{C\in\mathbf{\Pi}(A,B,D)\middle|\begin{vmatrix}B-C\\C-A\\E\end{vmatrix}\le 0\right\}.$$

Then, in the case $\varphi\in(0,\ \pi)$, the intersection of this plane with $\mathbf{S}(A,B,\varphi)$ is the union of two circular segments:

$$\check{\mathbf{S}}(A,B,D,\varphi)=\left\{C\in\check{\mathbf{\Pi}}(A,B,D)\backslash\{A,B\}\middle|\angle(B-C,C-A)\in[0,\varphi)\right\},$$

$$\hat{\mathbf{S}}(A,B,D,\varphi)=\left\{C\in\hat{\mathbf{\Pi}}(A,B,D)\backslash\{A,B\}\middle|\angle(B-C,C-A)\in[0,\varphi)\right\}.$$

Thus, in the case $\varphi\in(0,\ \pi)$, the set $\mathbf{S}(A,B,\varphi)$ is a three-dimensional body obtained by rotating either $\check{\mathbf{S}}(A,B,D,\varphi)$ or $\hat{\mathbf{S}}(A,B,D,\varphi)$ around the line $AB$, for any fixed $D\in\mathbb{R}^3\backslash AB$.

In this case:

- If $\varphi=\pi/2$, the set $\mathbf{S}(A,B,\pi/2\ )$ will be a sphere (see Fig. 1).
- If $\varphi\in(0,\pi/2)$, then $\mathbf{S}(A,B,\varphi)$ will be lemon-shaped (see Fig. 2).
- If $\varphi\in(\pi/2,\pi)$, then it will be apple-shaped (see Fig. 3).

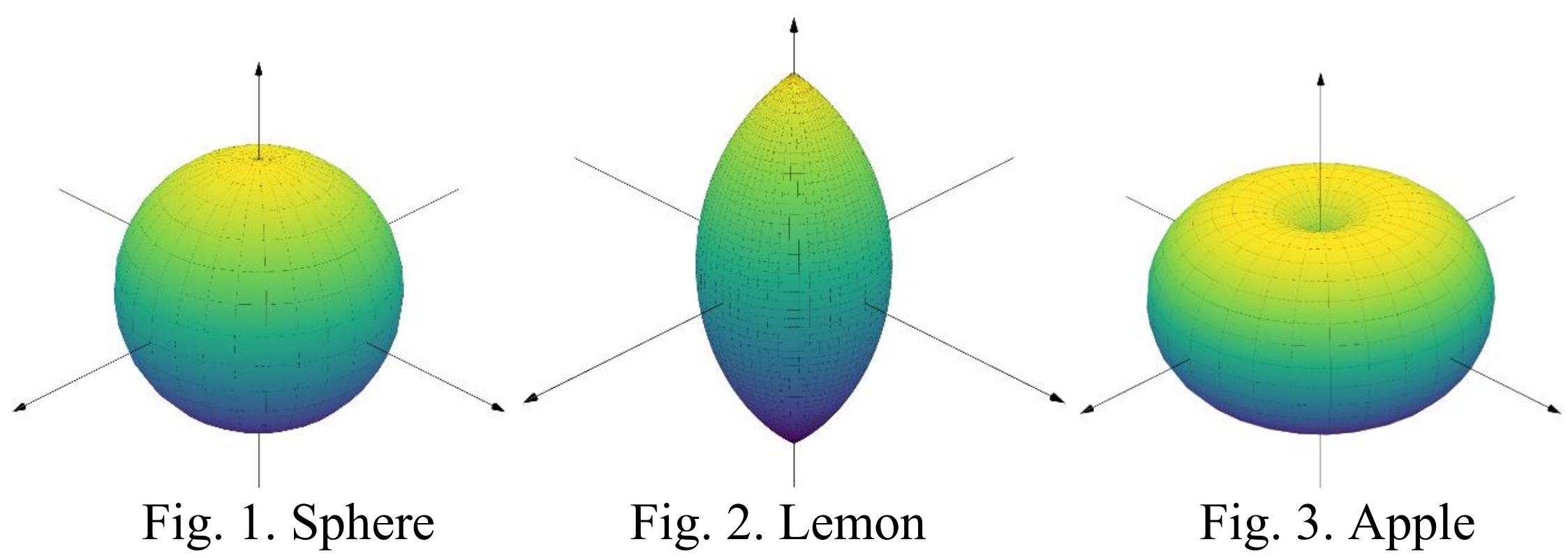

Fig. 1. Sphere　Fig. 2. Lemon　Fig. 3. Apple

If $\varphi>\pi$, then for convenience we set $\mathbf{S}(A,B,\varphi)=\mathbb{R}^3\backslash\{A,B\}$, which corresponds to equality (2). The case $\varphi=\pi$ is more subtle and is discussed in Remark 1.

**Remark 1.** Let us consider the case $\varphi=\pi$. To describe it, we introduce a parametrization of the line passing through the points $A$ and $B$: $A(t)=A+t(B-A)$, $t\in\mathbb{R}$. Then, by definition:

$$(A(-\infty),A]=\{A(t)|t\in(-\infty,0]\},\ (A(-\infty),A)=\{A(t)|t\in(-\infty,0)\},$$

$$[B,A(+\infty))=\{A(t)|t\in[1,+\infty)\},\ (B,A(+\infty))=\{A(t)|t\in(1,+\infty)\}.$$

It is directly verified [1] that

$$\mathbf{S}(A,B,\pi)=\mathbb{R}^3\backslash\{(A(-\infty),A]\cup[B,A(+\infty))\}\qquad(3)$$

(to justify (3), the reasoning for the two-dimensional case can be easily extended to the three-dimensional case we are considering).

Furthermore, for any non-zero vector $E\in\mathbb{R}^3$ and $\varphi\in(0,\pi]$, we define the angular sector:

$$\mathbf{C}(E,\varphi)=\{C\in\mathbb{R}^3|\angle(C,E)\in[0,\varphi)\},$$

assuming $\mathbf{0}\in\mathbf{C}(E,\varphi)$. We will also need the closures of the sets $\mathbf{S}(A,B,\varphi)$ and $\mathbf{C}(E,\varphi)$:

$$\operatorname{cl}\mathbf{S}(A,B,\varphi)=\{C\in\mathbb{R}^3|\angle(B-C,C-A)\in[0,\varphi]\},$$

$$\text{cl}\,\mathbf{C}(E,\varphi)=\{C\in\mathbb{R}^3 | \angle(C,E)\in[0,\varphi]\}.$$

## 2. Main Results

Let us consider a polygonal chain $\mathbf{L}$ connecting the points $A=B^{(0)}\in\mathbb{R}^3$ and $B=B^{(n+1)}\in\mathbb{R}^3$:

$$\mathbf{L}=\bigcup_{i=0}^{n}\left[B^{(i)},B^{(i+1)}\right], \qquad (4)$$

where $n\geq 1$. Let the turning angles at the vertices $B^{(i)}\in\mathbb{R}^3$ $(i=1,\dots,n)$ be bounded by a given number $\varphi\in(0,\pi)$ (in the case $n=1$, $\varphi=\pi$ is allowed):

$$\begin{gathered}\angle\left(B^{(i+1)}-B^{(i)},B^{(i)}-B^{(i-1)}\right)\in[0,\varphi], i=1,\dots,n,\\ B^{(i)}\neq B^{(i-1)}, i=0,1,\dots,n.\end{gathered} \qquad (5)$$

Condition (5) allows for the angle to be zero, which corresponds to motion without turning (the point $B^{(i)}$ effectively lies on the line connecting the adjacent segments). In such a situation, the actual number of turns becomes less than $n$. Since in this work we are interested in polygonal chains with exactly $n$ turns, the case with the additional condition is of greatest interest:

$$\angle\left(B^{(i+1)}-B^{(i)},B^{(i)}-B^{(i-1)}\right)\neq 0,\ i=1,\dots,n. \qquad (6)$$

In this regard, we note that in some problems there may be a restriction on the lengths of the segments. In such cases, zero turns (obtained as a result of splitting a segment that is too long) become permissible.

**Theorem 1.** *Let $\varphi>0$, $n\geq 1$, $n\varphi\leq\pi$, $A\neq B$, and let the polygonal chain (4) satisfy condition (5). Then:*

*1. For $n\geq 2$, we have $B^{(1)},\dots,B^{(n)}\in\mathbf{S}(A,B,n\varphi)$, and $B^{(i)}\neq B^{(j)}$ for any distinct $i,j\in\{0,1,\dots,n+1\}$.*

*2. For $n=1$, we have $B^{(1)}\in\text{cl}\,\mathbf{S}(A,B,\varphi)\setminus\{A,B\}$.*

The proof of Theorem 1 is given in Section 4.

**Lemma 1.** *Under the conditions of Theorem 1, a polygonal chain satisfying (4), (5) (and, in the case of interest to us, also condition (6)) exists:*

*– for any $B^{(1)}\in\mathbf{S}(A,B,n\varphi)$ when $n\geq 2$;*

*– for any $B^{(1)}\in\text{cl}\,\mathbf{S}(A,B,\varphi)\setminus\{A,B\}$ when $n=1$.*

The proof of Lemma 1 is given in Section 5.

**Remark 2.** Continuing Remark 1, we make clarifications for Lemma 1 and Theorem 1 in the case $n\varphi>\pi$. By analogy with the proof of Lemma 1, it can be shown that in this case, a polygonal chain $\mathbf{L}$ satisfying conditions (4), (5), (6) exists for any $B^{(1)}\in\mathbb{R}^3\setminus\{A,B\}$, which corresponds to the condition $\mathbf{S}(A,B,n\varphi)=\mathbb{R}^2\setminus\{A,B\}$ for this case. However, it should be noted that in this situation, the condition $B^{(i)}\neq B^{(j)}$ for any distinct $i,j\in\{0,1,\dots,n+1\}$ from Theorem 1 may not be satisfied.

**Remark 3.** As we see from Remark 2, the statements of Theorem 1 and Lemma 1 can be generalized to the case $n\varphi>\pi$ with appropriate clarifications (see Remark 2) and include the special case $n\varphi = \pi$. However, the most meaningful case is $n\varphi < \pi$, which is the focus of this work.

We denote by $\mathbf{S}^{(n)}(A, B, \varphi)$ the set of all possible sequences of interior vertices $\left(B^{(1)}, \dots, B^{(n)}\right)$ for which there exists a polygonal chain (4) satisfying condition (5). Using Theorem 1 and Lemma 1, we can provide the following description of the elements of the set $\mathbf{S}^{(n)}(A, B, \varphi)$.

**Theorem 2.** *Under the conditions of Theorem 1, we have:*

$$\left(B^{(1)}, \dots, B^{(n)}\right) \in \mathbf{S}^{(n)}(A, B, \varphi) \Leftrightarrow \left(B^{(1)}, \dots, B^{(n)}\right) \in$$
$$\in \mathbf{S}(A, B, n\varphi) \times \left\{\mathbf{S}\left(B^{(1)}, B, (n-1)\varphi\right) \cap \left[B^{(1)} + \mathrm{cl}\, \mathbf{C}\left(B^{(1)} - A, \varphi\right)\right]\right\} \times$$
$$\times \left\{\mathbf{S}\left(B^{(2)}, B, (n-2)\varphi\right) \cap \left[B^{(2)} + \mathrm{cl}\, \mathbf{C}\left(B^{(2)} - B^{(1)}, \varphi\right)\right]\right\} \times \cdots \times$$
$$\times \left\{\mathbf{S}\left(B^{(n-2)}, B, 2\varphi\right) \cap \left[B^{(n-2)} + \mathrm{cl}\, \mathbf{C}\left(B^{(n-2)} - B^{(n-3)}, \varphi\right)\right]\right\} \times$$
$$\times \left\{\left[\mathrm{cl}\, \mathbf{S}\left(B^{(n-1)}, B, \varphi\right) \setminus \{A, B\}\right] \cap \left[B^{(n-1)} + \mathrm{cl}\, \mathbf{C}\left(B^{(n-1)} - B^{(n-2)}, \varphi\right)\right]\right\}. \quad (7)$$

*Moreover, for any sequential selection of elements* $B^{(1)}, \dots, B^{(n)}$ *in accordance with the right-hand side of (7), the set from which the next element of this sequence is chosen will be non-empty for any set of preceding terms.*

The proof of Theorem 2 is given in Section 6. Formula (7) essentially uses the principle of dynamic programming [5–7].

Formula (7) can serve as a basis for constructing algorithms to enumerate admissible polygonal chains. If an additional constraint is imposed on all points – namely, that they must belong to some set $\mathbf{Q} \subset \mathbb{R}^3$ (for example, a coordinate parallelepiped or a uniform grid on it) – then an intersection with $\mathbf{Q}$ must be added to each term of the direct product on the right-hand side of (7). This makes it possible to formulate and solve optimization problems of the form:

$$f\left(B^{(1)}, \dots, B^{(n)}\right) \to \min; \left(B^{(1)}, \dots, B^{(n)}\right) \in \mathbf{S}^{(n)}(A, B, \varphi) \cap \mathbf{Q}^n \quad (8)$$

– exactly in the discrete case (when $\mathbf{Q}$ is finite), and approximately (after approximation) in the continuous case (see the corresponding constructions for the two-dimensional case in [1], Section 5).

## 3. Auxiliary concepts and facts

To prove the main statements, we will need some auxiliary propositions.

**Statement 2.** *Let* $A, B \in \mathbb{R}^3$, $A \neq B$, $\varphi \in (0, \pi)$. *Then*

$$\mathbf{S}(A, B, \varphi) \subset \mathrm{cl}\, \mathbf{S}(A, B, \varphi) \subset A + \mathbf{C}(B - A, \varphi).$$

**Proof.** Let $C \in \mathrm{cl}\, \mathbf{S}(A, B, \varphi)$. First, consider the case $C \in AB = \{A + t(B - A) | t \in \mathbb{R}\}$.

- If $C \in AB \setminus [A, B]$, then $\angle(B - C, C - A) = \pi$. Since $\varphi \in (0, \pi)$, it follows that $C \notin \mathrm{cl}\, \mathbf{S}(A, B, \varphi) = \{C \in \mathbb{R}^3 | \angle(B - C, C - A) \in [0, \varphi]\}$.
- Now let $C \in [A, B]$. Then $\angle(C - A, B - A) = 0$, i.e. $C - A \in \mathbf{C}(B - A, \varphi)$ and $C \in A + \mathbf{C}(B - A, \varphi)$.

Now let $C \notin AB$ (i.e., points $A, B, C$ are not collinear). Then

$$\angle(B - C, C - A) \in (0, \varphi].$$

Consider triangle $\triangle ABC$, as well as the angle

$$\Psi(C) = \angle(C - A, B - A) = \angle(B - A, C - A).$$

Note that

$$\Psi(C) < \angle(B - C, C - A) \le \varphi$$

(an exterior angle of a triangle is greater than any of its interior angles), whence

$$C - A \in \mathbf{C}(B - A, \varphi) \Rightarrow C \in A + \mathbf{C}(B - A, \varphi).\blacksquare$$

It is not difficult to show that for an angular sector in the considered three-dimensional case, many properties obtained for the two-dimensional case are preserved. The proof is similar, but now, instead of the Lemma on the additivity of angles, we use the triangle inequality for angles – i.e., Statement 1.

The following series of propositions describes the nesting property of such angular sectors in the sequential construction of a polygonal chain.

**Statement 3.** *Let* $A, B$ *be distinct points in* $\mathbb{R}^3$, *and* $\varphi \in (0, \pi)$. *Then*

$$B + \operatorname{cl} \mathbf{C}(B - A, \varphi) \subseteq A + \mathbf{C}(B - A, \varphi).$$

**Proof.** Let $C \in B + \operatorname{cl} \mathbf{C}(B - A, \varphi)$, i.e., $\alpha = \angle(C - B, B - A) \in [0, \varphi]$.

Consider two cases.

**Case 1:** $\alpha = 0$. Then $C - B = t(B - A)$ for some $t \ge 0$, which implies $C - A = (t + 1)(B - A)$ and $\angle(C - A, B - A) = 0 \in [0, \varphi)$.

**Case 2:** $\alpha \in (0, \varphi]$. In this case, the points $A, B, C$ are not collinear ($\alpha \notin \{0, \pi\}$) and form a triangle. Moreover,

$$\angle(C - A, B - A) < \alpha,$$

since the angle $\angle(C - A, B - A)$ is an interior angle, while the angle $\alpha$ is an exterior angle of the triangle $\triangle ABC$. Consequently, $\angle(C - A, B - A) \in (0, \varphi)$.

Thus, in both cases, $\angle(C - A, B - A) \in [0, \varphi)$, which means $C \in A + \mathbf{C}(B - A, \varphi)$.$\blacksquare$

**Remark 4.** Note that in the case $\varphi = \pi$, Statement 2 does not hold. However, in this case one can prove a weaker inclusion:

$$B + \mathbf{C}(B - A, \pi) \subset A + \mathbf{C}. \qquad (9)$$

Indeed, if $E \neq \mathbf{0}$, then

$$\forall C \in \mathbb{R}^3 \;\; \angle(C, E) = \pi \Leftrightarrow C \in (E(-\infty), 0) = \{-tE | t > 0\},$$

whence

$$\mathbf{C}(E, \pi) = \{C \in \mathbb{R}^3 | \angle(C, E) \in [0, \pi)\} = \mathbb{R}^3 \backslash (E(-\infty), 0).$$

But then, under the conditions of Statement 2, $\operatorname{cl} \mathbf{C}(B - A, \pi) = \mathbb{R}^3$, $B + \operatorname{cl} \mathbf{C}(B - A, \pi) = \mathbb{R}^3$, $A + \mathbf{C}(B - A, \pi) = \mathbb{R}^3 \backslash (A(-\infty), A) \neq \mathbb{R}^3$, where $(A(-\infty), A) = \{A + t(B - A) | t \in (-\infty, 0)\}$, i.e., the inclusion

$$B + \operatorname{cl} \mathbf{C}(B - A, \pi) \subseteq A + \mathbf{C}(B - A, \pi)$$

does not hold. At the same time,

$$B + \mathbf{C}(B - A, \pi) = \mathbb{R}^3 \backslash (A(-\infty), B), A + \mathbf{C}(B - A, \pi) = \mathbb{R}^3 \backslash (A(-\infty), A),$$
$$[B + \mathbf{C}(B - A, \pi)] \backslash [A + \mathbf{C}(B - A, \pi)] = \emptyset,$$
$$[A + \mathbf{C}(B - A, \pi)] \backslash [B + \mathbf{C}(B - A, \pi)] = [A, B),$$

i.e., (9) is proven. $\blacksquare$

**Statement 4.** *Let* $A, B$ *be distinct points in* $\mathbb{R}^3$, $\varphi_1, \varphi_2 > 0$, $\varphi_1 + \varphi_2 \in (0, \pi]$, $C \in B + \operatorname{cl} \mathbf{C}(B - A, \varphi_1)$, $C \neq B$ *and* $D \in C + \operatorname{cl} \mathbf{C}(C - B, \varphi_2)$. *Then*

$$D \in B + \mathbf{C}(B - A, \varphi_1 + \varphi_2) \subseteq A + \mathbf{C}(B - A, \varphi_1 + \varphi_2).$$

**Proof.** From the conditions of the statement to be proved, we obtain $\angle(C-B,B-A)\in[0,\varphi_1]$. By Statement 3, $D\in C+\mathrm{cl}\,\mathbf{C}(C-B,\varphi_2)\subseteq B+\mathbf{C}(C-B,\varphi_2)$, whence $\angle(D-B,C-B)\in[0,\varphi_2)$. Applying Statement 1 (since $C\neq B$), we get:

$$\angle(D-B,B-A)\leq\angle(D-B,C-B)+\angle(C-B,B-A)<\varphi_1+\varphi_2,$$

which means that $D\in B+\mathbf{C}(B-A,\varphi_1+\varphi_2)$. The second inclusion follows from Statement 3 as well as from Remark 4. ■

**Statement 5.** *Let* $A,B$ *be distinct points in* $\mathbb{R}^3$, $\varphi_1,\varphi_2,\varphi_3>0$, $\varphi_1+\varphi_2+\varphi_3\in(0,\pi]$, $C\in B+\mathrm{cl}\,\mathbf{C}(B-A,\varphi_1)$, $C\neq B$, $D\in C+\mathrm{cl}\,\mathbf{C}(C-B,\varphi_2)$, $D\neq C$, *and* $E\in D+\mathrm{cl}\,\mathbf{C}(D-C,\varphi_3)$. *Then*

$$E\in B+\mathbf{C}(B-A,\varphi_1+\varphi_2+\varphi_3)\subseteq A+\mathbf{C}(B-A,\varphi_1+\varphi_2+\varphi_3).$$

**Proof.** Applying Statement 4 to the points $B,C,D,E$ with angles $\varphi_2,\varphi_3$, we obtain $E\in C+\mathbf{C}(C-B,\varphi_2+\varphi_3)$. Then, applying Statement 4 to the points $A,B,C,E$ with angles $\varphi_1,\varphi_2+\varphi_3$, we get the desired result. ■

**Statement 6.** *Let* $m\geq 1$, $A=B^{(0)}$, $B^{(1)},\dots,B^{(m+1)}\in\mathbb{R}^2$, $B^{(i)}\neq B^{(i-1)}$ *for* $i=1,\dots,m+1$, *and for some* $\varphi_i>0$ *the following conditions hold:*

$$B^{(i+1)}\in B^{(i)}+\mathrm{cl}\,\mathbf{C}\big(B^{(i)}-B^{(i-1)},\varphi_i\big),\ i=1,\dots,m.$$

*Then:*

1. *If* $\sum_{i=1}^m\varphi_i\leq\pi$, *then the following inclusion holds:*
$$B^{(m+1)}\in B^{(1)}+\mathbf{C}\big(B^{(1)}-A,\textstyle\sum_{i=1}^m\varphi_i\big)\subseteq A+\mathbf{C}\big(B^{(1)}-A,\textstyle\sum_{i=1}^m\varphi_i\big).$$
2. *In the case* $0\leq j<i\leq m$ *and* $\varphi_{j+1}+\dots+\varphi_i\leq\pi$, *the following holds:*
$$B^{(i+1)}\in B^{(j+1)}+\mathbf{C}\big(B^{(j+1)}-B^{(j)},\varphi_{j+1}+\dots+\varphi_i\big)\subseteq$$
$$\subseteq B^{(j)}+\mathbf{C}\big(B^{(j+1)}-B^{(j)},\varphi_{j+1}+\dots+\varphi_i\big).$$

**Proof.** Statement 1) is easily proved by induction on $m$.

**Base case** ($m=1$) is obvious (see Statement 3).

**Inductive step:** apply the induction hypothesis to the points $B^{(1)},\dots,B^{(m+1)}$, obtaining

$$B^{(m+1)}\in B^{(2)}+\mathbf{C}\big(B^{(2)}-B^{(1)},\textstyle\sum_{i=2}^m\varphi_i\big).$$

Moreover, by the conditions, $B^{(2)}\in B^{(1)}+\mathrm{cl}\,\mathbf{C}\big(B^{(1)}-A,\varphi_1\big)$. Then apply Statement 4 to the points $A$, $B^{(1)}$, $B^{(2)}$, $B^{(m+1)}$ with angles $\varphi_1$ and $\sum_{i=2}^m\varphi_i$, which yields the required inclusion for $B^{(m+1)}$. Statement 2) is a consequence of Statement 1). ■

Corollary of Statement 6 is:

**Statement 7.** *Suppose we are under the conditions of Statement 6. Then, if* $\sum_{i=1}^m\varphi_i\leq\pi$, *it follows that* $B^{(m+1)}\neq B^{(0)}$. *Furthermore, in the case* $0\leq j<i\leq m$ *and* $\varphi_{j+1}+\dots+\varphi_i\in(0,\pi]$, *the following holds:* $i+1\geq j+2\geq 2$ *and* $B^{(i+1)}\neq B^{(j)}$.

**Proof.** Let us prove that $B^{(m+1)}\neq B^{(0)}$ (the proof that $B^{(i+1)}\neq B^{(j)}$ is analogous). Suppose that $B^{(m+1)}=B^{(0)}$. Then, by Statement 6, we have

$$B^{(m+1)}=B^{(0)}\in B^{(1)}+\mathbf{C}\big(B^{(1)}-B^{(0)},\varphi_1+\dots+\varphi_m\big),$$

whence

$$B^{(0)} - B^{(1)} \in \mathbf{C}\big(B^{(1)} - B^{(0)}, \varphi_1 + \dots + \varphi_m\big).$$

However, this contradicts the fact that

$$\angle\big( B^{(0)} - B^{(1)}, B^{(1)} - B^{(0)}\big) = \pi \geq \varphi_1 + \dots + \varphi_m. \blacksquare$$

**Remark 5.** Statement 8 can be strengthened by replacing the condition $B^{(m+1)} \neq B^{(0)}$ with the weaker one:

$$B^{(m+1)} \neq B^{(1)} + t\big(B^{(0)} - B^{(1)}\big), \forall t > 0.$$

Accordingly, for $0 \leq j < i \leq m$ and $\varphi_{j+1} + \dots + \varphi_i \in (0, \pi)$, the condition $B^{(i+1)} \neq B^{(j)}$, where $i + 1 \geq j + 2 \geq 2$, can be replaced with:

$$B^{(i+1)} \neq B^{(j+1)} + t\big(B^{(j)} - B^{(j+1)}\big), \forall t > 0.$$

## 4. Proof of Theorem 1

We will prove Theorem 1 by induction on $n$ – the number of turns.

**Base case.** Let $n = 1$. Condition (5) gives $\angle\big(B - B^{(1)}, B^{(1)} - A\big) \in [0, \varphi]$, which, by definition, means that $B^{(1)} \in \mathrm{cl}\,\mathbf{S}(A, B, \varphi)\backslash\{A, B\}$. The statement is proved.

**Inductive step.** Suppose that for some $n \geq 2$, the theorem holds in all cases with the number of turns $\leq n - 1$.

We will show that

$$B^{(1)} \in \mathbf{S}(A, B, n\varphi). \qquad (10)$$

If (10) does not hold, we have:

$$\Psi\big(B^{(1)}\big) = \angle\big(B - B^{(1)}, B^{(1)} - A\big) \geq n\varphi. \qquad (11)$$

By the induction hypothesis, the following holds:

$$B^{(2)}, \dots, B^{(n)} \in \mathbf{S}\big(B^{(1)}, B, (n-1)\varphi\big),$$

and in the case $n = 2$, $B^{(2)} \in \mathrm{cl}\,\mathbf{S}\big(B^{(1)}, B, \varphi\big)\backslash\{A, B\}$. In any of these cases, by Statement 2, we obtain:

$$B^{(2)} \in B^{(1)} + \mathbf{C}\big(B - B^{(1)}, (n-1)\varphi\big),$$

and consequently,

$$\angle\big(B^{(2)} - B^{(1)}, B - B^{(1)}\big) < (n-1)\varphi. \qquad (12)$$

On the other hand, by (5),

$$\angle\big(B^{(2)} - B^{(1)}, B^{(1)} - A\big) \leq \varphi. \qquad (13)$$

From (12) and (13), by the triangle inequality for angles (Statement 1) and the fact that, by (5), $B^{(2)} \neq B^{(1)}$, we get:

$$\angle\big(B - B^{(1)}, B^{(1)} - A\big) \leq \angle\big(B - B^{(1)}, B^{(2)} - B^{(1)}\big) +$$
$$+\angle\big(B^{(2)} - B^{(1)}, B^{(1)} - A\big) < n\varphi,$$

which contradicts (11).

We will now prove that

$$B^{(2)} \in \mathbf{S}(A, B, n\varphi)$$

i.e.,

$$\angle\big(B - B^{(2)}, B^{(2)} - A\big) < n\varphi.$$

Suppose that

$$\angle\big(B - B^{(2)}, B^{(2)} - A\big) \geq n\varphi. \qquad (14)$$

By the induction hypothesis, we have:

$$B^{(3)}, \dots, B^{(n)} \in \mathbf{S}\big(B^{(2)}, B, (n-2)\varphi\big),$$

and in the case $n = 3$, $B^{(3)} \in \operatorname{cl} \mathbf{S}\big(B^{(2)}, B, \varphi\big)\backslash\{A, B\}$. In any of these cases, by Statement 2, we obtain:

$$B^{(3)} \in B^{(2)} + \mathbf{C}\big(B - B^{(2)}, (n-2)\varphi\big),$$

and consequently,

$$\angle\big(B^{(3)} - B^{(2)}, B - B^{(2)}\big) < (n-2)\varphi. \qquad (15)$$

On the other hand, by (5),

$$\angle\big(B^{(3)} - B^{(2)}, B^{(2)} - B^{(1)}\big) \leq \varphi,\ \angle\big(B^{(2)} - B^{(1)}, B^{(1)} - A\big) \leq \varphi. \quad (16)$$

From (15) and (16), using (5) and by the triangle inequality for angles (Statement 1), we get:

$$\angle\big(B - B^{(2)}, B^{(2)} - A\big) \leq \angle\big(B - B^{(2)}, B^{(3)} - B^{(2)}\big) +$$
$$+\angle\big(B^{(3)} - B^{(2)}, B^{(2)} - B^{(1)}\big) + \angle\big(B^{(2)} - B^{(1)}, B^{(1)} - A\big) < n\varphi,$$

which contradicts (14).

Similarly, we prove that $B^{(3)}, \dots, B^{(n)} \in \mathbf{S}(A, B, n\varphi)$.

Now note that, by (5),

$$B^{(i+1)} \in B^{(i)} + \operatorname{cl} \mathbf{C}\big(B^{(i)} - B^{(i-1)}, \varphi\big),\ i = 1, \dots, n,$$

and at the same time $n\varphi \leq \pi$, i.e., for $m = n, \varphi_i = \varphi, i = 1, \dots, m$, we are under the conditions of Statements 6 and 7, which imply that $B^{(i)} \neq B^{(j)}$ for any distinct $i, j \in \{0,1, \dots, n+1\}$ ($B^{(i)} \neq B^{(i-1)}, i = 0,1, \dots, n,$ by (5); and if $0 \leq j < i \leq n$, then $B^{(i+1)} \neq B^{(j)}$ by Statement 7). ■

**Remark 5.** Suppose we are under the conditions of Theorem 1 and, moreover, $n\varphi = \pi$. Then, by this theorem, it follows that $B^{(1)}, \dots, B^{(n)} \in \mathbf{S}(A, B, \pi)$, whence, by (3), we obtain:

$$B^{(i)} \in \mathbb{R}^3\backslash\{(A(-\infty), A] \cup [B, A(\infty))\}, i = 1, \dots, n.$$

## 5. Proof of Lemma 1

We will now present the **proof** of Lemma 1. Consider the non-trivial case when $n \geq 2$. Let $D = B^{(1)} \in \mathbf{S}(A, B, n\varphi)$, where $n\varphi \leq \pi$.

**Case 1.** Let $D = B^{(1)} \notin AB = \{A + t(B - A) | t \in \mathbb{R}\}$. In this case, the points $A, B, D$ are not collinear. Let $\mathbf{\Pi}(A, B, D)$ be the plane passing through the points $A, B, D$, i.e. (see [4]):

$$\mathbf{\Pi}(A, B, D) = \operatorname{aff}\{A, B, D\} = \{t_1 A + t_2 B + t_3 D | t_1, t_2, t_3 \in \mathbb{R}, t_1 + t_2 + t_3 = 1\}.$$

The line $AB$ divides the plane $\mathbf{\Pi}(A, B, D)$ into two half-planes: $\check{\mathbf{\Pi}}(A, B, D)$ – conventionally the right one, and $\hat{\mathbf{\Pi}}(A, B, D)$ – conventionally the left one, which share a common boundary – the line $AB$. For example, if $E \in \mathbb{R}^3$ is a non-zero vector perpendicular to the plane $\mathbf{\Pi}(A, B, D)$, then we can take $\check{\mathbf{\Pi}}(A, B, D)$ as the set of all points $C \in \mathbf{\Pi}(A, B, D)$ such that either $C \in AB$, or the vectors $B - C, C - A$ and $E$ form a right-handed triple of vectors, i.e.:

$$\begin{vmatrix} B-C \\ C-A \\ E \end{vmatrix} > 0.$$

This leads to the following formulas:

$$\check{\mathbf{\Pi}}(A,B,D) = \left\{ C \in \mathbf{\Pi}(A,B,D) \middle| \begin{vmatrix} B-C \\ C-A \\ E \end{vmatrix} \geq 0 \right\},$$

$$\hat{\mathbf{\Pi}}(A,B,D) = \left\{ C \in \mathbf{\Pi}(A,B,D) \middle| \begin{vmatrix} B-C \\ C-A \\ E \end{vmatrix} \leq 0 \right\}.$$

Then, in the case $n\varphi < \pi$, the intersection of this plane with $\mathbf{S}(A,B,n\varphi)$ is the union of two circular segments:

$$\check{\mathbf{S}}(A,B,D,n\varphi) = \{C \in \check{\mathbf{\Pi}}(A,B,D)\backslash\{A,B\} | \angle(B-C,C-A) \in [0,n\varphi)\},$$

$$\hat{\mathbf{S}}(A,B,D,n\varphi) = \{C \in \hat{\mathbf{\Pi}}(A,B,D)\backslash\{A,B\} | \angle(B-C,C-A) \in [0,n\varphi)\}.$$

In this case, $D = B^{(1)} \in \check{\mathbf{S}}(A,B,D,n\varphi)$ or $D = B^{(1)} \in \hat{\mathbf{S}}(A,B,D,n\varphi)$. Both of these cases are treated similarly and, by a simple linear orthogonal transformation of variables, are reduced to the two-dimensional case considered in [1]. Since the transition matrix under such a transformation is orthogonal, the values of angles between vectors are preserved. Therefore, the solution for the planar case described in [1], when transformed back to the three-dimensional case, will yield the desired solution.

Let now $n\varphi = \pi$. Then, since $B^{(1)} \in \mathbf{S}(A,B,\pi)$, we have

$$\alpha = \angle\left(B - B^{(1)}, B^{(1)} - A\right) \in (0,\pi),\ \bar{\alpha} = \tfrac{\alpha+\pi}{2} \in (\alpha,\pi),\ B^{(1)} \in \mathbf{S}(A,B,\bar{\alpha}),$$

and now the case under consideration is reduced to the one considered above by replacing the expression $n\varphi$ with $\bar{\alpha}$ in $\mathbf{S}(A,B,n\varphi)$, $\check{\mathbf{S}}(A,B,D,n\varphi)$, and $\hat{\mathbf{S}}(A,B,D,n\varphi)$.

**Case 2.** Now consider the case when $B^{(1)} \in AB$, $n\varphi \leq \pi$. Then:

- In the case $n\varphi < \pi$, it is easily shown in [1] that $B^{(1)} \in (A,B)$.
- In the case $n\varphi = \pi$, this follows from equality (3).

Now let $D$ be an arbitrary point from $\mathbb{R}^3\backslash AB$. Then, as in Case 1, using the plane $\mathbf{\Pi}(A,B,D)$, we reduce the problem to considering the planar case described in [1] (corresponding to Case 1 in [1]). The solution obtained for the planar case, when transformed back to the three-dimensional setting, will yield the desired solution.

## 6. Proof of Theorem 2

**Proof.** Recall that by $\mathbf{S}^{(n)}(A,B,\varphi)$ we denote the set of all possible sequences of internal vertices $\left(B^{(1)}, \dots, B^{(n)}\right)$ for which there exists a polygonal chain $\mathbf{L}$ satisfying (4), (5).

**Direct inclusion.** Let $\left(B^{(1)}, \dots, B^{(n)}\right) \in \mathbf{S}^{(n)}(A,B,\varphi)$, i.e. $B^{(1)}, \dots, B^{(n)}$ is a sequence of internal points of the polygonal chain $\mathbf{L}$ satisfying (4), (5). Consider the non-trivial case with $n \geq 2$. Applying Theorem 1 to the polygonal chain **L**, we obtain $B^{(1)} \in \mathbf{S}(A,B,n\varphi)$. Next, applying the same theorem to the part of the polygonal chain starting at $B^{(1)}$, we have $B^{(2)} \in \mathbf{S}\big(B^{(1)}, B, (n-1)\varphi\big)$. Moreover, from condition (5) for $i = 1$ it follows that $\angle\big(B^{(2)} - B^{(1)}, B^{(1)} - A\big) \in [-\varphi, \varphi]$, i.e. $B^{(2)} \in B^{(1)} + \mathrm{cl}\,\mathbf{C}\big(B^{(1)} - A, \varphi\big)$. Thus,

$$B^{(2)} \in \mathbf{S}\big(B^{(1)}, B, (n-1)\varphi\big) \cap \big[B^{(1)} + \mathrm{cl}\,\mathbf{C}\big(B^{(1)} - A, \varphi\big)\big].$$

Continuing in the same manner, for each $i = 2, \dots, n-2$ we obtain

$$B^{(i+1)} \in \mathbf{S}\big(B^{(i)}, B, (n-i)\varphi\big) \cap \big[B^{(i)} + \mathrm{cl}\,\mathbf{C}\big(B^{(i)} - B^{(i-1)}, \varphi\big)\big],$$

and for the last point,

$$B^{(n)} \in \big[\mathrm{cl}\,\mathbf{S}\big(B^{(n-1)}, B, \varphi\big) \backslash \{A, B\}\big] \cap \big[B^{(n-1)} + \mathrm{cl}\,\mathbf{C}\big(B^{(n-1)} - B^{(n-2)}, \varphi\big)\big].$$

Therefore, the sequence $\big(B^{(1)}, \dots, B^{(n)}\big)$ belongs to the direct product on the right-hand side of (7).

**Reverse inclusion.** We will prove that for any sequential selection of points $B^{(1)}, \dots, B^{(n)}$ in accordance with the right-hand side of (7), each member in this direct product will be non-empty and any such sequence satisfies (5).

Indeed, suppose we have chosen an arbitrary point

$$B^{(1)} \in \mathbf{S}(A, B, n\varphi). \tag{17}$$

By Lemma 1, for the point $B^{(1)}$ there exists a polygonal chain

$$\boldsymbol{L} = \big[B^{(0)}, B^{(1)}\big] \cup \big[B^{(1)}, B^{(2)}\big] \cup \dots \cup \big[B^{(n)}, B^{(n+1)}\big],$$

satisfying (5). Then, by (5) and by the definition of the set $\mathrm{cl}\,\mathbf{C}(E, \varphi)$, for the point $B^{(2)}$ from any such polygonal chain it holds that

$$B^{(2)} \in B^{(1)} + \mathrm{cl}\,\mathbf{C}\big(B^{(1)} - A, \varphi\big).$$

Note that the point $B^{(2)}$ is the first turning point of the polygonal chain

$$\bar{\boldsymbol{L}} = \big[B^{(1)}, B^{(2)}\big] \cup \dots \cup \big[B^{(n)}, B^{(n+1)}\big],$$

and, by Theorem 1,

$$B^{(2)} \in \mathbf{S}\big(B^{(1)}, B, (n-1)\varphi\big).$$

Thus,

$$B^{(2)} \in \mathbf{S}\big(B^{(1)}, B, (n-1)\varphi\big) \cap \big[B^{(1)} + \mathrm{cl}\,\mathbf{C}\big(B^{(1)} - A, \varphi\big)\big], \tag{18}$$

and the set on the right-hand side of (18) is non-empty.

Suppose we have already chosen two points $B^{(1)}, B^{(2)}$ satisfying (17), (18). By Lemma 1, for the point $B^{(2)}$ there exists a polygonal chain

$$\bar{\mathbf{L}} = \big[B^{(1)}, B^{(2)}\big] \cup \dots \cup \big[B^{(n)}, B^{(n+1)}\big],$$

satisfying (5). Then, by (5) and by the definition of the set $\mathrm{cl}\,\mathbf{C}(E, \varphi)$, for the point $B^{(3)}$ from any such polygonal chain it holds that

$$B^{(3)} \in B^{(2)} + \mathrm{cl}\,\mathbf{C}\big(B^{(2)} - B^{(1)}, \varphi\big).$$

Note that the point $B^{(3)}$ is the first turning point of the polygonal chain

$$\bar{\bar{\mathbf{L}}} = \big[B^{(2)}, B^{(3)}\big] \cup \dots \cup \big[B^{(n)}, B^{(n+1)}\big],$$

and, by Theorem 1,

$$B^{(3)} \in \mathbf{S}\big(B^{(2)}, B, (n-2)\varphi\big).$$

Thus,

$$B^{(3)} \in \mathbf{S}\big(B^{(2)}, B, (n-2)\varphi\big) \cap \big[B^{(2)} + \mathrm{cl}\,\mathbf{C}\big(B^{(2)} - B^{(1)}, \varphi\big)\big], \tag{19}$$

and the set on the right-hand side of (19) is non-empty.

In exactly the same way, it can be proved that each subsequent member in the direct product (7) is non-empty for any selection of the preceding members according to (7). Thus, the sequence of points $B^{(1)}, \dots, B^{(n)}$ in the direct product (7) can always be extended for any choice of the initial members (where $i \in \{0, \dots, n-1\}$). The satisfaction of (5) for this sequence follows from the condition

$$B^{(i)} \in B^{(i-1)} + \mathrm{cl}\,\mathbf{C}\big(B^{(i-1)} - B^{(i-2)}, \varphi\big),\ i = 2, \dots, n+1. \blacksquare$$

## Conclusion

In this work, two key statements have been proven – Theorem 1 and Lemma 1 – which describe the set of admissible piecewise linear polygonal chains connecting two given points in three-dimensional Euclidean space, subject to constraints on the number of turns and the magnitudes of angles at the breakpoints.
Based on these results, an explicit description (Theorem 2) of the set of all possible sequences of angular points of a polygonal chain has been obtained.
The resulting structural formula (7) can serve as a foundation for constructing algorithms to enumerate admissible polygonal chains. In the discrete case, when the points of the polygonal chain must belong to a given finite set (for instance, an integer lattice), this approach allows obtaining an exact solution to an optimization problem of the form (8).
In the continuous case, the formula opens up the possibility of constructing finite approximations of the set of admissible polygonal chains with a prescribed accuracy, and subsequently finding approximate solutions to optimization problems.
Issues of set approximation and optimization methods are discussed in detail in [1, 2, 8].
Future research may be directed towards developing specific algorithmic implementations of the proposed approach.
If necessary, a further generalization of the obtained results to an arbitrary finite-dimensional Euclidean space is possible. With such a generalization, all the reasoning remains valid.